\documentclass[12pt]{article}
\usepackage{amsmath,amsthm,amsfonts}
\usepackage{url}

\newcommand{\pf}{\noindent {\bf Proof:} }
\newtheorem{theorem}{Theorem}
\newtheorem{remark}{Remark}

\newtheorem{corollary}[theorem]{Corollary}
\newtheorem{lemma}[theorem]{Lemma}

\newtheorem{definition}[theorem]{Definition}
\title{Dejean's conjecture holds for $n\ge 30$}
\author{James Currie\thanks{The author is
supported by an NSERC Discovery Grant.} and
Narad Rampersad\thanks{The author is supported by an NSERC Postdoctoral
Fellowship.} \\
Department of Mathematics and Statistics \\
University of Winnipeg \\
515 Portage Avenue \\
Winnipeg, Manitoba R3B 2E9 (Canada) \\
\url{j.currie@uwinnipeg.ca} \\
\url{n.rampersad@uwinnipeg.ca}}

\begin{document}
\date{\today}
\maketitle

\begin{abstract}
We extend Carpi's results by showing that Dejean's conjecture
holds for $n\ge 30.$
\end{abstract}

\section{Introduction}

Repetitions in words have been studied starting with Thue \cite{thueI,thue} at the beginning of the previous century. Much study has also been given to repetitions with fractional exponent \cite{brandenburg,carpi,dejean,longfrac,krieger,mignosi}. If $n>1$ is an integer, then an {\bf $n$-power} is a non-empty word $x^n$, i.e., word $x$ repeated $n$ times in a row. For rational $r>1$, a {\bf fractional $r$-power} is a non-empty word $w=x^{\lfloor r\rfloor}x'$ such that $x'$ is the prefix of $x$ of length $(r-\lfloor r\rfloor)|x|$. For example, $01010$ is a $5/2$-power. A basic problem is that of identifying the {\it repetitive threshold} for each alphabet size $n>1$: \begin{quote}What is the infimum of $r$ such that an infinite sequence on $n$ letters exists, not containing any $r$-powers?\end{quote}
We call this infimum the {\it repetitive threshold} of an $n$-letter alphabet, denoted by $RT(n)$. Dejean's conjecture \cite{dejean} is that
$$RT(n)=\left\{\begin{array}{ll}7/4,&n=3\\
7/5,&n=4\\
n/(n-1)&n\ne 3,4\end{array}\right.$$

The values $RT(2)$, $RT(3)$, $RT(4)$ were established by Thue, Dejean and Pansiot, respectively \cite{thueI,dejean,pansiot}. Moulin-Ollagnier \cite{ollagnier} verified Dejean's conjecture for $5\le n\le 11$, while Mohammad-Noori and Currie \cite{morteza} proved the conjecture for $12\le n\le 14$.

An exciting new development has recently occurred with the work of Carpi \cite{carpi}, who showed that Dejean's conjecture holds for $n\ge 33$. Verification of the conjecture is now only lacking for a finite number of values. In the present paper, we sharpen Carpi's methods to show that Dejean's conjecture holds for $n\ge 30$.

\section{Preliminaries}

The following definitions are from sections 8 and 9 of
\cite{carpi}: Fix $n\ge 30$. Let $m=\lfloor(n-3)/6\rfloor$. Let
$A_m =\{1,2,\ldots, m\}$. Let ker $\psi=\{v\in A_m^*|\forall a\in
A_m$, 4 {\it divides }$|v|_a\}.$ (In fact, this is not Carpi's
{\it definition} of ker $\psi$, but rather the assertion of his
Lemma~9.1.) A word $v\in A_m^+$ is a $\psi$-{\bf kernel
repetition} if it has period $q$ and a prefix $v'$ of length $q$
such that $v'\in$ ker $\psi$, $(n-1)(|v|+1)\ge nq-3.$

It will be convenient to have the following new definition: If $v$
has period $q$ and its prefix $v'$ of length $q$ is in ker $\psi$,
we say that $q$ is a {\bf kernel period} of $v$.

As Carpi states at the beginning of section 9 of \cite{carpi}:
\begin{quotation}
By the results of the previous sections, at least in the case
$n\ge 30$, in order to construct an infinite word on $n$ letters
avoiding factors of any exponent larger than $n/(n-1)$, it is
sufficient to find an infinite word on the alphabet $A_m$ avoiding
$\psi$-kernel repetitions.
\end{quotation}
For $m=5$, Carpi produces such an infinite word, based on a
paper-folding construction. He thus establishes Dejean's
conjecture for $n\ge 33$. In the present paper, we give an
infinite word on the alphabet $A_4$ avoiding $\psi$-kernel
repetitions. We thus establish Dejean's conjecture for $n\ge 30$.

\begin{definition}Let $f:A_4^*\rightarrow A_4^*$ be
defined by $f(1)=121$, $f(2)=123$, $f(3)=141$, $f(4)=142$. Let $w$
be the fixed point of $f$.
\end{definition}
It is useful to note that the frequency matrix of $f$, i.e., 
$$[|f(i)|_j]_{4\times 4}=\left[\begin{array}{cccc}2&1&0&0\\1&1&1&0\\2&0&0&1\\1&1&0&1\end{array}\right]$$ 
 has an
inverse modulo 4.
\begin{remark}\label{remark}
Let $q$ be a non-negative integer, $q\le 1966$. Fix $n=32$.
\begin{enumerate}\item[R1:]{Word $w$ contains no $\psi$-kernel repetition $v$ with kernel period $q$.}
\item[R2:]{Word $w$ contains no factor $v$ with
 kernel period $q$ such that $|v|/q\ge 35/34$.}
\end{enumerate}
\end{remark}

Note that ${32\over 31}-{34\over 31q}={35\over 34}$ when
$q={34^2\over 3}=385{1\over 3}$, so neither piece of the remark
implies the other. Note also that the conditions of the remark
become {\bf less} stringent for $n=30,31$. One also verifies that
$${35\over 34}+{9\over 2(1967)}\le{32\over 31}-{34\over 31q}$$ for $q\ge 1967$.
To show that $w$ contains no $\psi$-kernel repetitions for $n=$
30, 31, 32, it thus suffices to verify R1 and to show that word
$w$ contains no factor $v$ with kernel period $q\ge 1967$ such
that
\begin{equation}\label{stronger}|v|/q\ge 35/34 +
9/2(1967).\end{equation}

The remarks R1 and R2 are verified by computer search, so we will
consider the second part of this attack. Fix $q\ge 1967$, and
suppose that $v$ is a factor of $w$ with kernel period $q$, and
$|v|/q\ge 35/34$. Without loss of generality, suppose that no
extension of $v$ has period $q$. Write $v=sf(u)p$ where $s$ (resp.
$p$) is a suffix (resp. prefix) of the image of a letter, and
$|s|$ $(\mbox{ resp. }|p|)$ $\le 2$.

If $|v|\le q+2$, then $35/34\le (q+2)/q$ and $1/34\le 2/q$,
forcing $q\le 68$. This contradicts R2. We will therefore assume
that $|v|\ge q+3$.

Suppose $|s|=2$. Since $|v|\ge q+3$, write $v=s0zs0v'$, where
$|s0z|=q$. Examining $f$, we see that the letter $a_s$ preceding
any occurrence of $s0$ in $w$ is uniquely determined if $|s|=2$.
It follows that $a_sv$ is a factor of $w$ with kernel period $q$,
contradicting the maximality of $v$. We conclude that $|s|\le 1$.

Again considering $f$, we see that if $t$ is any factor of $w$ of
length 3, and $u_1t$, $u_2t$ are prefixes of $w$, then
$|u_1|\equiv |u_2|$ (mod 3). Since $|v|\ge q+3$, we conclude that
3 divides $q$. Write $q=3q_0$. Since $|s|\le 1$, $|p|\le 2$ and
$|v|\ge q+3$, we see that $|f(u)|\ge q.$ Thus $f(u)$ has a prefix
of length $q=3q_0$ which is in ker $\psi$. As the frequency matrix
of $f$ is invertible modulo 4, the prefix of $u$ of length $q_0$
is in ker $\psi$. We see that
$${|v|\over q}\le{3|u|+3\over 3q_0}={|u|\over q_0}+{1\over q_0}.$$
\begin{lemma}
Let $s$ be a non-negative integer. If factor $v$ of $w$ has kernel
period $q$, where $q\le 1966(3^s)$, then
$${|v|\over q}< {35\over 34}+{3\over
1966}\sum_{j=0}^{s-1}3^{-j}.$$
\end{lemma}
\pf If $s=0$, this is implied by R2. Suppose $t>0$ and the result
holds for $s<t$. Suppose that $1966(3^{t-1})<q\le 1966(3^{t})$ and
there is a factor $v$ of $w$ such that $v$ has kernel period $q$.
Suppose that $|v|/q\ge 35/34$. Without loss of generality, suppose
that no extension of $v$ has period $q$. We have seen that there
is a factor $u$ of $w$ with kernel period $q_0=q/3$,
$1966(3^{t-2})<q_0\le 1966(3^{t-1})$ such that
\begin{eqnarray*}
|v|/q&\le&|u|/q_0+1/q_0\\
&<&\left({35\over 34}+{3\over
1966}\sum_{j=0}^{t-2}3^{-j}\right)+{1\over
q_0}\mbox{ (by the induction hypothesis)}\\
&<&{35\over 34}+{3\over 1966}\sum_{j=0}^{t-2}3^{-j}+{1\over
1966(3^{t-2})}\\
 &=&{35\over 34}+{3\over
1966}\sum_{j=0}^{t-2}3^{-j}+{3\over 1966(3^{t-1})}\\
 &=&{35\over 34}+{3\over 1966}\sum_{j=0}^{t-1}3^{-j}.\Box\end{eqnarray*}
\begin{theorem}
Word $w$ contains no factor $v$ with kernel period $q$ such that
$$|v|/q\ge 35/34 + 9/2(1966).$$
\end{theorem}
\pf Suppose that factor $v$ of $w$ has kernel period $q$ such that
(\ref{stronger}) holds. By Remark~\ref{remark}, we have $q\ge
1966$. By the previous lemma, for some non-negative $s$,
$$|v|/q<{35\over 34}+{3\over 1966}\sum_{j=0}^{s-1}3^{-j}<{35\over
34}+{3\over 1966}\sum_{j=0}^\infty 3^{-j}={35\over 34}+{9\over
2(1966)}.\Box$$
\begin{corollary}Dejean's conjecture holds for $n=30, 31, 32$.
\end{corollary}
\section*{Appendix: Computer search}
Suppose that some factor $v$ of $w$ has kernel period $q\le
1966$ and either $31(|v|+1)\ge 32q-3$ or $|v|/q\ge
35/34+9/2(1967).$ Without loss of generality, taking such a $v$ as
short as possible, we may assume that
$$|v|\le\left\lceil {32(1966)-3\over 31}-1\right\rceil=2029.$$
$$\mbox{(We also have }
\left\lceil (1966)\left({35\over 34}+{9\over
2(1967)}\right)\right\rceil=2029.)$$ If $|v|>3$, $v$ is a factor
of $f(u)$ for some factor $u$ of $w$ with $|u|\le (|v|+4)/3.$ For a
non-negative integer $r$, let $g(r)=\lfloor (r+4)/3\rfloor$. Since
$g^7(2029)=2<3$, (here the exponent denotes iterated function
composition) word $v$ must be a factor of $f^7(u)$ for some factor
$u$ of $w$, $|u|=2$.

The word $u_0=23141121142$ contains all 8 factors of $w$
which have length 2. To establish R1 and R2, one thus checks that
they hold for the single word $f^7(u_0)$ (which is of length 24,057).
In fact, we performed this computer check on the word $f^7(u_1)$,
where $u_1=11421231211231411$ contains all 13 factors of $w$
which have length 3.

\end{document}